\documentclass[11pt]{article}
\usepackage{enumerate}
\usepackage{amsmath}
\usepackage{amssymb}
\usepackage{amsthm}
\usepackage{amscd}
\usepackage{catmac}
\DeclareSymbolFontAlphabet{\Bbb}{AMSb}
\newcommand{\id}{\text{{\rm id}}}

\newcommand{\norm}[1]{ \|#1 \| }
\newcommand{\bignorm}[1]{\left \|#1 \right \| }
\newcommand{\convex}[1]{\mathbf{M^{(#1)}}}
\newcommand{\concave}[1]{\mathbf{M_{(#1)}}}
\newcommand{\cotype}[1]{\mathbf{C_{#1}}}
\newcommand{\type}[1]{\mathbf{T_{#1}}}

\newcommand{\K}{\Bbb{K}}
\newcommand{\C}{\Bbb{C}}
\newcommand{\R}{\Bbb{R}}

\newcommand{\LL}{{\cal L}}

\newcommand{\Id}{\hookrightarrow}

\newcommand{\tensor}[1]{\tilde{\otimes}_{#1}}
\newcommand{\ui}{\mathcal{S}}
\theoremstyle{definition}
\newtheorem{defin}{Definition}

\theoremstyle{plain}
\newtheorem{lemma}[defin]{Lemma}

\newtheorem{cor}[defin]{Corollary}
\newtheorem{prop}[defin]{Proposition}
\newtheorem{Approx}[defin]{Approximation Lemma}
\theoremstyle{remark}
\numberwithin{equation}{section}
\begin{document}
\bibliographystyle{amsalpha}
\title{\bf A complex interpolation formula for tensor products of 
vector-valued Banach function spaces}
\author{Andreas Defant and Carsten Michels 
\thanks{The second named author was supported by a GradF\"oG Stipendium
 of the Land Niedersachsen (State of Lower Saxony). Mathematics 
Subject Classification (1991): Primary 46M35; Secondary 46M05, 46E40, 46B70.}
}
\date{}
\maketitle
\begin{abstract}
We prove the complex interpolation formula
$$
[X_0(E_0) \tensor{\varepsilon} Y_0(F_0), X_1(E_1) \tensor{\varepsilon} 
Y_1(F_1)]_\theta = [X_0(E_0),X_1(E_1)]_\theta \tensor{\varepsilon} 
[Y_0(F_0),Y_1(F_1)]_\theta, 
$$
for the injective tensor product of vector-valued Banach function
 spaces $X_i(E_i)$ and $Y_i(F_i)$ satisfying certain geometric assumptions. 
This result unifies 
results of Kouba, and moreover, our approach offers an alternate proof of 
 Kouba's interpolation formula for scalar-valued Banach function spaces.
\end{abstract}
The following theorem for the complex 
interpolation of injective tensor products of vector-valued Banach 
function spaces is proved:
\\[10pt]
{\bf Theorem.} {\em
Let $X_0(\mu),X_1(\mu),Y_0(\nu),Y_1(\nu)$ be  real-valued Banach function 
spaces, and  $[E_0,E_1]$ and $[F_0,F_1]$ interpolation couples 
of complex Banach spaces with dense intersections. Then  
for  \mbox{$0 < \theta<1$} the equality
\begin{equation}
\label{vke}
[X_0(E_0) \tensor{\varepsilon} Y_0(F_0), X_1(E_1) \tensor{\varepsilon} 
Y_1(F_1)]_\theta = [X_0(E_0),X_1(E_1)]_\theta \tensor{\varepsilon} 
[Y_0(F_0),Y_1(F_1)]_\theta 
\end{equation}
holds algebraically and topologically 
whenever the Banach lattices $X_0,X_1,Y_0,Y_1$ are $2$-concave and the 
Banach spaces $E_i$ and $F_i$ satisfy one of the following conditions:
\\[5pt]
$(1)$ $E_0', E_1',F_0'$ and $F_1'$ are type~$2$ spaces.
\\
$(2)$ $E_0',E_1'$ are type~$2$ spaces and $F_0=F_1$ is a cotype~$2$ space.
\\
$(3)$ $E_0=E_1$ and $F_0=F_1$ are cotype~$2$ spaces.  
}
\\[10pt]
This is an extension and unification of deep results due to Kouba 
\cite{kouba} who proved
 the preceding interpolation formula if one of the couples 
$[X_0,X_1]$ and $[E_0,E_1]$, and one of the couples $[Y_0,Y_1]$ and 
$[F_0,F_1]$ is trivial (i.\,e. either $X_0=X_1=\R$ or 
$E_0=E_1=\C$, and either $Y_0=Y_1=\R$ or $F_0=F_1=\C$). Moreover, 
following an 
idea of Pisier \cite{pisier90} and based on variants of the 
Maurey--Rosenthal Factorization Theorem (see \cite{defant}), our approach
 offers an alternate proof of Kouba's interpolation formula for 
complex-valued Banach function spaces: For $2$-concave complex-valued 
Banach function spaces $X_0(\mu),X_1(\mu),Y_0(\nu),Y_1(\nu)$ and $0<\theta<1$
\begin{equation}
\label{oldkouba}
[X_0 \tensor{\varepsilon} Y_0,X_1 \tensor{\varepsilon} Y_1]_\theta 
= [X_0,X_1]_\theta \tensor{\varepsilon} [Y_0,Y_1]_\theta.
\end{equation}
\\[5pt]
The main ingredients of the proof  will 
be ``uniform estimates'' of 
\begin{equation}
\label{firstdef}
d_\theta[M_0,M_1]:=\norm{\LL(\ell_2,[M_0,M_1]_\theta) \Id 
[\LL(\ell_2,M_0),\LL(\ell_2,M_1)]_\theta},
\end{equation}
where $[M_0,M_1]$ is an interpolation couple of two 
$n$-dimensional Banach spaces. Such 
estimates proved to be of independent interest: The facts 
$\sup_n d_\theta[\ell_1^n,\ell_2^n]< \infty $ 
(see \cite{pisier90} and \cite[3.5]{kouba}; here it is a consequence of 
Proposition~\ref{semikouba}) and its non-commutative analogue
$\sup_n d_\theta[\ui_1^n,\ui_2^n] < \infty$ for finite-dimensional 
Schatten classes
 (due to Junge in \cite[4.2.6]{junge} and based on an extension of Kouba's
 formulas for the Haagerup tensor product of operator spaces due to
 \cite{pisier96}) were used in \cite{dm98} in order to study
 so-called ``Bennett--Carl Inequalities'' for identity operators
 between finite-dimensional symmetric Banach sequence spaces 
 as well as their ``non-commutative analogues'' 
for identity operators between finite-dimensional unitary ideals.  
\section*{Preliminaries}
 We shall use standard notation and
 notions from Banach space theory, as presented e.\,g. in 
\cite{djt}, \cite{lt77}, \cite{lt} and \cite{tj}; for tensor products 
of Banach spaces we refer to \cite{df}.
If $E$ is a Banach space, then
 $B_E$ is its (closed) unit ball and $E'$ its dual, and $FIN(E)$ 
stands for the collection of all its finite-dimensional subspaces. 
As usual $\LL(E,F)$
 denotes the Banach space of all (bounded and linear) operators from
 $E$ into $F$ endowed with the operator norm $\norm{\cdot}$.
For a Banach space $E$ of type~$2$ (resp. cotype~$2$) we write 
$\type{2}(E)$ (resp. $\cotype{2}(E)$) for its (Rademacher) type~$2$ constant
 (resp. cotype~$2$ constant), and for $1 \le r \le \infty$ we denote
 by $\convex{r}(X)$ (resp. $\concave{r}(X)$) the $r$-convexity 
(resp. $r$-concavity) constant of an $r$-convex (resp. $r$-concave) 
Banach lattice $X$. Recall that for Banach spaces $E,F$  the injective 
 norm on $E \otimes F$ is defined by
$$
\norm{z}_{E \otimes_\varepsilon F}:=\sup \{|\langle x' \otimes y',z 
\rangle | \, | \, x'\in B_{E'}, y' \in B_{F'} \}, \quad z \in E \otimes F,
$$
and with $E \tensor{\varepsilon} F$ we denote the completion of 
$E \otimes F$ endowed with this norm. We will extensively use the fact that 
the equality $E \otimes_\varepsilon F= \LL(E',F)$ holds isometrically 
whenever one of the two involved spaces is finite-dimensional.
\par
Let $(\Omega,\Sigma,\mu)$ be a $\sigma$-finite and complete measure space,
 and denote all ($\mu$-a.e. equivalence classes of) real-valued measurable
 functions on $\Omega$ by $L_0(\mu)$. A Banach space $X=X(\mu)$ of 
  functions in $L_0(\mu)$ is said to be a 
Banach function space if it satisfies the following conditions:
\\[5pt]
\begin{samepage}
(I) If $|f| \le |g|$, with $f \in L_0(\mu)$ and $g \in X(\mu)$, then
 $f \in X(\mu)$ and $\norm{f}_X \le \norm{g}_X$.
\\[2pt]
(II) For every $A \in \Sigma$ with $\mu(A)<\infty$ the characteristic 
function $\chi_A$ of $A$ belongs to $X(\mu)$.
\end{samepage}
\\[10pt]
A finite-dimensional real Banach space $X=(\R^n,\norm{\cdot}_X)$ is called 
an 
$n$-dimensional lattice if  $\norm{\cdot}_X$
 is a lattice norm; clearly, $X$ then is a Banach function space in the 
above sense. 
For Banach function spaces $X_0(\mu),X_1(\mu)$  and $0<\theta<1$ 
the space $X_0^{1-\theta} X_1^\theta$ is defined to be the set of 
 functions $f \in L_0(\mu)$ for which there exist
 $g \in X_0$ and $h \in X_1$ such that $|f|=|g|^{1-\theta} \cdot
 |h|^\theta$. Together with the norm 
$$
\norm{f}_{X_0^{1-\theta} X_1^\theta}:= \inf \{ 
\norm{g}_{X_0}^{1-\theta} \cdot \norm{h}_{X_1}^\theta \, | \, 
|f|=|g|^{1-\theta} \cdot |h|^\theta, g \in X_0, h \in X_1 \},
$$
$X_0^{1-\theta} X_1^\theta$ becomes a Banach function space (with 
respect to $(\Omega,\Sigma,\mu)$). It can be easily seen (see e.\,g. 
\cite[p.~218/219]{tj}) that if for $1 \le r < \infty$ the lattices 
$X_0$ and $X_1$ are both $r$-convex or both $r$-concave, then 
$X_0^{1-\theta} X_1^\theta$ also has this property, with
\begin{equation}
\label{convexinterpol} 
\convex{r}(X_0^{1-\theta}X_1^\theta) \le \convex{r}(X_0)^{1-\theta} 
\cdot \convex{r}(X_1)^\theta,
\end{equation} 
\begin{equation}
\label{concaveinterpol} 
\concave{r}(X_0^{1-\theta}X_1^\theta) \le \concave{r}(X_0)^{1-\theta} 
\cdot \concave{r}(X_1)^\theta,
\end{equation}
respectively.  
\\[10pt]
Let $X(\mu)$ be a Banach function space and $E$ a Banach space. A function
 $x$ defined on $\Omega$ with values in $E$ is said to be strongly 
measurable if there exists a sequence of strictly simple functions on 
$\Omega$ converging to $x$ almost everywhere; here a function $y$ on 
$\Omega$ with values in $E$ is called strictly simple if it assumes 
only finitely many non-zero values, each on a measurable set with finite
 measure. Then by $X(E)$ we denote the collection of all strongly measurable
 functions $x$ with values in $E$ for which $\norm{x(\cdot)}_E \in X$. 
Together with the norm $\norm{x}_{X(E)}:=\norm{\norm{x(\cdot)}_E}_X$, 
this vector space becomes a Banach space ($\K$-linear whenever $E$ is 
$\K$-linear). 
\par 
For all information on complex interpolation we refer to \cite{BL} and 
\cite{kps}. Given a (complex) interpolation couple $[E_0,E_1]$, we write
 $E_\Delta:=E_0 \cap E_1$, and as usual denote for $0 < \theta <1$ 
 the complex interpolation space with respect to 
$[E_0,E_1]$ and $\theta$ by 
$[E_0,E_1]_\theta$. If we speak of a finite-dimensional interpolation
 couple $[E_0,E_1]$, this always means that both spaces have the same
 finite dimension. Clearly, if $[E_0,E_1]$ is an interpolation 
couple and $X_0(\mu),X_1(\mu)$ are Banach function spaces, then 
$[X_0(E_0),X_1(E_1)]$ is an interpolation couple. We will heavily use the 
following complex interpolation formula due to Calder\'{o}n 
\cite[13.6]{calderon}: For $0 < \theta <1$
\begin{equation}
\label{calderonformula}
[X_0(E_0),X_1(E_1)]_\theta = 
(X_0^{1-\theta} X_1^\theta)([E_0,E_1]_\theta)
\end{equation}
holds isometrically whenever $X_0$ or $X_1$ is $\sigma$-order continuous; 
note 
that under the assumptions of the theorem all involved Banach function 
spaces are $\sigma$-order continuous (for an argument see Section~4), and 
clearly this is true for finite-dimensional lattices. 
\section{The approximation lemma}
First we show---similar to \cite[Section~4]{kouba}---that  
equalities as stated in the above theorem are of a finite-dimensional
nature.
In order to make the following more readable, let us introduce the following
 notation: If $[E_0,E_1]$ is an interpolation couple, $E \subset E_\Delta$
 a subspace which is dense in $E_0,E_1$ and $\mathcal{A} \subset FIN(E)$ is cofinal
 (i.\,e. for every $G \in FIN(E)$ there exists $M \in \mathcal{A}$ with $G \subset M$),
 then the triple $([E_0,E_1],E,\mathcal{A})$ is called a {\em cofinal interpolation
 triple}. For $M \in FIN(E)$ we denote by $M_0$ (resp. $M_1$) the
 subspace $M$ of $E_0$ (resp. $E_1$) endowed with the induced norm. 
\par
The following two lemmas are crucial. The first one is an only slight 
modification of \cite[4.1]{kouba}; we omit its proof. 
\begin{lemma}
\label{lifting}
Let $([E_0,E_1],E,\mathcal{A})$ be a cofinal interpolation triple and
$0<\theta<1$.  Then
 for each $\varepsilon>0$ and  $G \in FIN(E)$ there exists 
$M \in \mathcal{A}$ such that $G \subset M$ and for all $x \in G$
\begin{equation}
\label{liftinge}
 (1-\varepsilon)\cdot 
\norm{x}_{[M_0,M_1]_\theta} 
\le \norm{x}_{[E_0,E_1]_\theta} \le 
\norm{x}_{[M_0,M_1]_\theta}.
\end{equation}
\end{lemma}
 If $[M_0,M_1]$ and $[N_0,N_1]$ are finite-dimensional interpolation 
couples, then we define for $0 <\theta<1$
$$
d_\theta[M_0,M_1;N_0,N_1]:=\norm{
[M_0,M_1]_\theta \otimes_\varepsilon [N_0,N_1]_\theta \Id 
[M_0 \otimes_\varepsilon N_0,M_1 \otimes_\varepsilon N_1]_\theta}. 
$$ 
The second lemma---which for obvious reasons is called 
``approximation lemma''---reduces the proof of Kouba type formulas 
\eqref{vke} or \eqref{oldkouba} to uniform estimates of 
$d_\theta[M_0,M_1;N_0,N_1]$ for cofinally many suitable finite-dimensional 
subspaces of the underlying infinite-dimensional spaces. Its proof
 is very close to the proof of \cite[4.2]{kouba}, but we state it 
for the convenience of the reader.
\begin{Approx}
\label{epscor}
Let $([E_0,E_1],E,\mathcal{A})$ and $([F_0,F_1],F,\mathcal{B})$ be cofinal interpolation 
triples and $0<\theta<1$. If 
$$
d_\theta[E_0,E_1;F_0,F_1] :=\sup_{M\in \mathcal{A}} \sup_{N \in \mathcal{B}} 
d_\theta[M_0,M_1;N_0,N_1] < \infty,
$$
then 
$$
[E_0 \tensor{\varepsilon} F_0, E_1 \tensor{\varepsilon} F_1]_\theta 
= [E_0,E_1]_\theta \tensor{\varepsilon} [F_0,F_1]_\theta.
$$
\end{Approx}
\proof
 From the density assumptions we conclude that $E \otimes F$ is dense in
 $[E_0,E_1]_\theta \tensor{\varepsilon} [F_0,F_1]_\theta$ and in
 $[E_0 \tensor{\varepsilon} F_0,E_1 \tensor{\varepsilon} F_1]_\theta$,
 hence it is sufficient to show that for a given $z \in E \otimes F$
\begin{align}
\label{toshow1}
\norm{z}_{[E_0,E_1]_\theta \tensor{\varepsilon} [F_0,F_1]_\theta} & \le 
\norm{z}_{[E_0 \tensor{\varepsilon} F_0,E_1 \tensor{\varepsilon} 
F_1]_\theta} \\
\label{toshow2}
&\le d_\theta[E_0,E_1;F_0,F_1] \cdot 
\norm{z}_{[E_0,E_1]_\theta \tensor{\varepsilon} [F_0,F_1]_\theta}.
\end{align}
We start with a simple observation to show \eqref{toshow1}.
If $[M_0,M_1]$ and $[N_0,N_1]$ are finite-dimensional interpolation couples,
then  
\begin{equation}
\label{reverse}
\norm{[\LL(M_0,N_0),\LL(M_1,N_1)]_\theta \Id \LL([M_0,M_1]_\theta,
[N_0,N_1]_\theta)} \le 1;
\end{equation}
indeed, consider for $i=0,1$ the bilinear mapping
$$
\phi_i : \LL(M_i,N_i) \times M_i \rightarrow N_i, \quad (T,x) \mapsto Tx,
$$
which clearly has norm $1$, hence \eqref{reverse} follows from the fact 
that by bilinear interpolation (see \cite[4.4.1]{BL}) the interpolated mapping 
$$
\phi_\theta: [\LL(M_0,N_0), \LL(M_1,N_1)]_\theta \times [M_0,M_1]_\theta 
\rightarrow [N_0,N_1]_\theta
$$
also has norm $\le 1$. Now \eqref{toshow1} is a straightforward consequence:
Obviously $C:=\{M\otimes N \, 
| \, M \in \mathcal{A}, N \in \mathcal{B} \} \subset FIN(E \otimes F)$ is cofinal, hence, 
by Lemma~\ref{lifting} and the fact that the injective norm respects 
subspaces, there exist $M \in \mathcal{A}$ and $N \in \mathcal{B}$ such that 
$z \in M \otimes N$ and
$$
\norm{z}_{[M_0 \otimes_\varepsilon N_0,M_1 \otimes_\varepsilon N_1]_\theta}
 \le (1+\varepsilon) \cdot \norm{z}_{[E_0 \tensor{\varepsilon} F_0,E_1 
\tensor{\varepsilon} F_1]_\theta}.
$$
Finally, by the mapping property of the injective norm and \eqref{reverse},
\begin{align*}
\norm{z}_{[E_0,E_1]_\theta \tensor{\varepsilon} [F_0,F_1]_\theta} &\le 
\norm{z}_{[M_0,M_1]_\theta \otimes_\varepsilon [N_0,N_1]_\theta}
 \\ 
&\le \norm{z}_{[M_0 \otimes_\varepsilon N_0,M_1 
\otimes_\varepsilon N_1]_\theta} \\
& \le (1+\varepsilon) \cdot \norm{z}_{[E_0 \tensor{\varepsilon} 
F_0,E_1 \tensor{\varepsilon} F_1]_\theta}.
\end{align*}
In order to show \eqref{toshow2} let 
$z \in G \otimes H$ for some $G \in FIN(E), H \in FIN(F)$, and
 choose by Lemma~\ref{lifting} subspaces $M \in \mathcal{A}$ and 
$N \in \mathcal{B}$ such that $G \subset M, H \subset N$ and
$$
\norm{(G,\norm{\cdot}_{[E_0,E_1]_\theta}) \Id [M_0,M_1]_\theta} \le 
\sqrt{1+\varepsilon},
$$
$$
\norm{(H,\norm{\cdot}_{[F_0,F_1]_\theta}) \Id [N_0,N_1]_\theta} \le 
\sqrt{1+\varepsilon}.
$$
Then, by the mapping property,
$$
\norm{(G,\norm{\cdot}_{[E_0,E_1]_\theta}) \otimes_\varepsilon 
 (H,\norm{\cdot}_{[F_0,F_1]_\theta}) \Id [M_0,M_1]_\theta 
\otimes_\varepsilon [N_0,N_1]_\theta} \le 1+\varepsilon,
$$
hence, since the injective norm respects subspaces,
$$
\norm{z}_{[M_0,M_1]_\theta \otimes_\varepsilon [N_0,N_1]_\theta} 
 \le (1 +\varepsilon) \cdot
 \norm{z}_{[E_0,E_1]_\theta \otimes_\varepsilon [F_0,F_1]_\theta}.
$$
By the usual interpolation theorem  we obtain
\begin{align*}
\norm{z}_{[E_0 \tensor{\varepsilon} F_0,E_1 \tensor{\varepsilon} 
F_1]_\theta} & \le \norm{z}_{[M_0 \otimes_\varepsilon N_0, M_1 
\otimes_\varepsilon N_1]_\theta} \\
&\le d_\theta [M_0,M_1;N_0,N_1] \cdot
\norm{z}_{[M_0,M_1]_\theta \otimes_\varepsilon [N_0,N_1]_\theta} \\
& \le  (1 +\varepsilon) \cdot d_\theta [E_0,E_1;F_0,F_1] \cdot 
 \norm{z}_{[E_0,E_1]_\theta \otimes_\varepsilon [F_0,F_1]_\theta}. 
\end{align*}
\qed 
\section{The Hilbert space case}
 Recall for a finite-dimensional interpolation couple $[E_0,E_1]$ the
 definition of $d_\theta[E_0,E_1]$ from \eqref{firstdef}, 
and  note that by the approximation lemma
$$
d_\theta[E_0,E_1]=\sup_n d_\theta[\ell_2^n,\ell_2^n;E_0,E_1].
$$
The main step in the proof of \eqref{vke}  is  the following estimate:
\begin{prop}
\label{semikouba}
Let $X_0,X_1$ be  $n$-dimensional lattices and $[E_0,E_1]$ 
 a finite-dimensional interpolation couple. Then for each $0 < \theta<1$
\begin{equation}
\label{semivkouba}
\begin{split}
d_\theta[X_0(E_0) & ,X_1(E_1)] \\
&\le \sqrt{2} \cdot \cotype{2}([E_0,E_1]_\theta) 
\cdot \concave{2}(X_0)^{1-\theta} \cdot \concave{2}(X_1)^\theta  \cdot 
d_\theta[\ell_2^n(E_0),\ell_2^n(E_1)].
\end{split}
\end{equation}
\end{prop}
Before giving the proof we collect some  facts about 
so-called  powers of 
finite-dimensional lattices. For $0<r<\infty$
 and an $n$-dimensional lattice $X$ with $\convex{\max(1,r)}(X)=1$ 
(recall that $\convex{1}(X)=1$) 
$$
\norm{x}_r := \norm{|x|^{1/r}}_X^r, \quad x \in \R^n
$$
defines a lattice norm on $\R^n$ (see e.\,g. \cite{defant}); 
the $n$-dimensional lattice $(\R^n,\norm{\cdot}_r)$ will be 
denoted by $X^r$.
\begin{lemma}
\label{observations}
Let $X,X_0,X_1$ be  $n$-dimensional lattices, $E$  
a Banach space, $\lambda \in \R^n$ and \mbox{$0<\theta<1$}.
\vspace{-10pt} 
\begin{enumerate}[(a)]
\item
\label{mgtensorid}
If $\concave{2}(X)=1$, then
$\norm{D_\lambda \otimes \id: \ell_2^n(E) \rightarrow X(E)} \le 
\norm{D_\lambda} =
 \norm{\lambda}_{(((X')^2)')^{1/2}}$, where $D_\lambda: \ell_2^n \rightarrow X$
 denotes the diagonal operator associated with $\lambda$.
\item
$(X_0^{1-\theta} X_1^\theta)'= (X_0')^{1-\theta} (X_1')^{\theta}$
 isometrically.
\item
For $0<r<\infty$ let $\convex{\max(1,r)}(X_0)=\convex{\max(1,r)}(X_1)=1$. 
Then 
$\left (X_0^{1-\theta} X_1^\theta \right)^r=(X_0^r)^{1-\theta} 
 (X_1^r)^\theta$ 
 isometrically.
\end{enumerate}
\end{lemma}
\proof (a)  For $x \in \ell_2^n(E)$ \vspace{-5pt}
$$
\norm{(D_\lambda \otimes \id)x}_{X(E)} = \norm{(\lambda_k \cdot \norm{x_k})_k}_X
 \le \norm{D_\lambda: \ell_2^n \rightarrow X} \cdot
\left(\sum_{k=1}^n \norm{x_k}^2 \right)^{1/2},
$$
\vspace{-5pt}
and (note that $\convex{2}(X')=\concave{2}(X)=1$)
\begin{align*}
\norm{\lambda}_{(((X')^2)')^{1/2}} & = \norm{\lambda^2}_{((X')^2)'}^{1/2} =
 \sup_{\norm{\mu}_{(X')^2} \le 1} \norm{\lambda^2  \mu}_{\ell_1^n}^{1/2} 
 = \sup_{\norm{|\mu|^{1/2}}_{X'} \le 1} \norm{\lambda^2  
\mu}_{\ell_1^n}^{1/2}\\
& = \sup_{\norm{\nu}_{X'} \le 1} \norm{\lambda  \nu}_{\ell_2^n} 
 = \sup_{\norm{\nu}_{X'} \le 1} \sup_{\norm{\mu}_{\ell_2^n}\le 1} \left |
 \sum_{i=1}^n \lambda_i \nu_i  \mu_i \right | \\ 
&= \sup_{\norm{\mu}_{\ell_2^n}\le 1} \sup_{\norm{\nu}_{X'} \le 1}  \left |
 \sum_{i=1}^n \lambda_i \mu_i  \nu_i  \right | 
 = \sup_{\norm{\mu}_{\ell_2^n} \le 1} \norm{\lambda  \mu }_X = 
\norm{D_\lambda: \ell_2^n \rightarrow X}. 
\end{align*}
\\[5pt]
(b) 
By the Calder\'{o}n formula \eqref{calderonformula}, 
the duality theorem \cite[4.5.2]{BL} and the 
fact that $Y(\C)'=Y'(\C)$ holds isometrically for every finite-dimensional 
lattice
 $Y$, one arrives at the isometric identity
$$
(X_0^{1-\theta} X_1^\theta)'(\C)=((X_0')^{1-\theta} 
(X_1')^\theta)(\C),
$$
 which clearly implies the above statement.
\\[5pt] (c) First note that  $\convex{\max(1,r)}(X_0^{1-\theta} 
X_1^\theta)=1$ by \eqref{convexinterpol}, hence the power  
$(X_0^{1-\theta} X_1^\theta)^r$ is normed.
Let $V:=(X_0^{1-\theta} X_1^\theta)^r$ and 
$W:= (X_0^r)^{1-\theta}  (X_1^r)^\theta$. Then, if $|f|^{1/r}
=|g|^{1-\theta} \cdot |h|^\theta$,
$$
\norm{f}_W  \le \norm{|g|^r}_{X_0^r}^{1-\theta} \cdot 
\norm{|h|^r}_{X_1^r}^\theta = \left ( \norm{g}_{X_0}^{1-\theta} \cdot
 \norm{h}_{X_1}^\theta \right)^r,
$$
which clearly implies $\norm{f}_W \le \norm{f}_V$. Conversely, let 
$|f|=|g|^{1-\theta} \cdot |h|^\theta$. Then
\begin{samepage}
$$
\norm{f}_V = \norm{|f|^{1/r}}_{X_0^{1-\theta} X_1^\theta}^r 
\le \norm{|g|^{1/r}}_{X_0}^{r(1-\theta)} \cdot 
\norm{|h|^{1/r}}_{X_1}^{r \theta}= \norm{g}_{X_0^r}^{1-\theta} \cdot
 \norm{h}_{X_1^r}^\theta,
$$
hence $\norm{f}_V \le \norm{f}_W$. \qed
\end{samepage}
\par     
Another important tool for the proof of 
\eqref{semivkouba} is a  variant of the 
Maurey--Rosenthal Factorization Theorem (\cite{maurey}) 
 for vector-valued Banach function spaces given in \cite{defant}. 
\begin{lemma}
\label{factorization}
Let $X(\mu)$ be a $2$-concave Banach function space and $E$ a Banach 
space of cotype~$2$. Then each $T \in \LL(\ell_2,X(E))$ factorizes as 
follows:
\begin{center}
\resetparms
\Vtriangle[\ell_2`X(E)`L_2(\mu,E);T`R`M_g \otimes \id],
\end{center}
where $R: \ell_2 \rightarrow L_2(\mu,E)$ is an operator and 
$M_g: L_2(\mu) \rightarrow X$ a multiplication operator with respect to
$g \in L_0(\mu)$ such that  
$\norm{R} \cdot \norm{M_g} \le \sqrt{2} \cdot \cotype{2}(E)  
 \cdot \concave{2}(X) \cdot \norm{T}$.  
\end{lemma}
\proof
Let $D_n:=\{-1,+1\}^n$, $\mu_n(\{\omega\}):=1/{2^n}$ for $\omega \in D_n$ and
 $\varepsilon_i: D_n \rightarrow \{-1,+1\}$ the $i$-th canonical projection.
 Then for $x_1, \ldots,x_n \in \ell_2$
\begin{equation*}
\label{maureyineq}
\begin{split}
\bignorm{ \left ( \sum_{i=1}^n \norm{Tx_i(\cdot)}_E^2 \right)^{1/2} }_X
 & \le \sqrt{2} \cdot \cotype{2}(E) \cdot \bignorm{ \int_{D_n} \norm{
 \sum_{i=1}^n \varepsilon_i(\omega) \cdot Tx_i (\cdot)}_E \, d\mu_n(\omega)}_X
 \\
& \le \sqrt{2} \cdot \cotype{2}(E) \cdot \int_{D_n} \left \| 
\norm{ (\sum_{i=1}^n \varepsilon_i(\omega) \cdot Tx_i ) (\cdot)}_E \right
 \|_X \, d\mu_n(\omega) \\
&= \sqrt{2} \cdot \cotype{2}(E) \cdot \int_{D_n} \bignorm{ T ( \sum_{i=1}^n
 \varepsilon_i(\omega) \cdot x_i)}_{X(E)} \, d\mu_n(\omega) \\
&\le \sqrt{2} \cdot \cotype{2}(E) \cdot  \norm{T} \cdot
\left ( \sum_{i=1}^n \norm{x_i}_{\ell_2}^2 \right )^{1/2}
\end{split}
\end{equation*}
(the constant $\sqrt{2}$ comes from the Khinchine--Kahane inequality for
 the case ``$L_2$ versus $L_1$''),
hence by \cite[4.4]{defant} there exists  $0 \le \omega \in L_0(\mu)$
 with
\begin{equation}
\label{weight2}
\sup_{y \in B_{L_2(\mu)}} \norm{\omega^{1/2} \cdot y}_X \le 
\sqrt{2} \cdot \concave{2}(X) \cdot \cotype{2}(E) \cdot \norm{T}
\end{equation}
such that for all $x \in \ell_2$
\begin{equation}
\label{ie}
\left ( \int_\Omega \norm{Tx(\cdot)}_E^2 /\omega \, d\mu \right)^{1/2}
 \le \norm{x}_{\ell_2}. 
\end{equation}
Define the operator $R \in \LL(\ell_2,L_2(\mu,E))$
 by $Rx:= Tx/{\omega^{1/2}}$ for $x \in \ell_2$ (well-defined by \eqref{ie}) 
and the multiplication operator $M_g: L_2(\mu) 
\rightarrow X$ with $g:=\omega^{1/2}$ (well-defined by \eqref{weight2}). 
Clearly, this produces the desired factorization. \qed
\par
Now we are prepared for the {\em Proof} of Proposition~\ref{semikouba}. 
Its main idea---the use of factorizations of Maurey--Rosenthal type---is 
taken from \cite{pisier90}.
\\[10pt] 
 Without loss of generality we may assume that $\concave{2}(X_0)=
\concave{2}(X_1)=1$; 
indeed, let $Y_0$ and $Y_1$ be the associated renormed
 lattices such that $\concave{2}(Y_i)=1$ and
\mbox{$\norm{X_i \Id Y_i} \cdot \norm{Y_i \Id X_i} \le 
\concave{2}(X_i)$} for $i=0,1$ (see e.\,g. \cite[1.d.8]{lt}).
 Now consider the factorization
\begin{center}
\setsqparms[1`1`1`1;2000`500]
\square[\ell_2 \otimes_\varepsilon [X_0(E_0),X_1(E_1)]_\theta`
[\ell_2 \otimes_\varepsilon X_0(E_0),\ell_2 \otimes_\varepsilon 
X_1(E_1)]_\theta`
\ell_2 \otimes_\varepsilon [Y_0(E_0),Y_1(E_1)]_\theta`
[\ell_2 \otimes_\varepsilon Y_0(E_0),\ell_2 \otimes_\varepsilon 
Y_1(E_1)\symbol{93}_\theta;
\id \otimes \id`u:= \id \otimes \id`v:= \id \otimes \id`\id \otimes \id] 
\end{center}
and observe that $\norm{u}\cdot \norm{v} \le \concave{2}(X_0)^{1-\theta}
 \cdot \concave{2}(X_1)^\theta.$
\\[10pt]
Put $X_\theta:=X_0^{1-\theta} X_1^\theta$.
Since $[X_0(E_0),X_1(E_1)]_\theta = X_\theta([E_0,E_1]_\theta)$ holds
 isometrically
 (see \eqref{calderonformula}) and $\concave{2}(X_\theta)=1$ 
(see \eqref{concaveinterpol}), by Lemma~\ref{factorization} every operator 
$T \in \LL(\ell_2,X_\theta([E_0,E_1]_\theta))$ factors
\begin{center}
\resetparms
\Vtriangle[\ell_2`X_\theta([E_0,E_1]_\theta)`\ell_2^n([E_0,E_1]_\theta);
T`R`D_\lambda \otimes \id],
\end{center}
with $\norm{R} \cdot \norm{D_\lambda} \le \sqrt{2} \cdot \cotype{2}
([E_0,E_1]_\theta)
 \cdot \norm{T: \ell_2 \rightarrow [X_0(E_0),X_1(E_1)]_\theta}$. 
Define $Y_\eta := (((X_\eta')^2)')^{1/2}$ for 
$\eta=0,1,\theta$; by Lemma~\ref{observations}~(b),(c) and 
the Calder\'{o}n formula \eqref{calderonformula} we have 
$[Y_0(\C),Y_1(\C)]_\theta =  Y_\theta(\C)$.  
By Lemma~\ref{observations}~(a) the mapping 
$$
\Phi_\eta: Y_\eta(\C) \rightarrow \LL(\ell_2^n(E_\eta),X_\eta(E_\eta)), 
\quad \mu \mapsto  D_\mu \otimes \id
$$
has norm $\le 1$, and consequently the interpolated mapping
$$
[\Phi_0,\Phi_1]_\theta : [Y_0(\C),Y_1(\C)]_\theta \rightarrow 
V:=[\LL(\ell_2^n(E_0), X_0(E_0)), \LL(\ell_2^n(E_1),X_1(E_1))]_\theta
$$
has norm $\le 1$.  Moreover,
 by bilinear interpolation (\cite[4.4.1]{BL}) the mapping
$$
U \times V \rightarrow W, \quad (u,v) \mapsto v \circ u,
$$
where
$
U:=[\LL(\ell_2,\ell_2^n(E_0), \LL(\ell_2,\ell_2^n(E_1))]_\theta$  
 and  
$W:= [\LL(\ell_2,X_0(E_0)),\LL(\ell_2,X_1(E_1))]_\theta 
$, also has norm $\le 1$. Since by definition 
$\norm{R}_U \le d_\theta[\ell_2^n(E_0),\ell_2^n(E_1)] \cdot \norm{R}$, 
 we obtain altogether
\begin{samepage}
\begin{align*}
\norm{T}_W = \norm{(D_\lambda \otimes \id) \circ R}_W &\le \norm{R}_U 
\cdot \norm{D_\lambda \otimes \id}_V = \norm{R}_U \cdot 
\norm{[\Phi_0,\Phi_1]_\theta (\lambda)}_V \\
& \le d_\theta[\ell_2^n(E_0),\ell_2^n(E_1)] \cdot \norm{R}
 \cdot \norm{\lambda}_{Y_\theta} \\
& \le d_\theta[\ell_2^n(E_0),\ell_2^n(E_1)] \cdot \sqrt{2} \cdot 
\cotype{2}([E_0,E_1]_\theta) \cdot \norm{T},
\end{align*}
the desired inequality. \qed
\end{samepage}
\begin{samepage}
\par
A quick look at \eqref{semivkouba} reveals that in the case $E=E_0=E_1$ 
 one has
\begin{cor}
Let $X_0,X_1$ be $n$-dimensional lattices and 
$E$ a finite-dimensional
 normed space. Then for $0<\theta<1$
\begin{equation}
\label{equal}
d_\theta[X_0(E),X_1(E)] \le \sqrt{2} \cdot \cotype{2}(E) 
\cdot \concave{2}(X_0)^{1-\theta} \cdot \concave{2}(X_1)^\theta. 
\end{equation}
\end{cor}
\end{samepage}
For the case that $E_0$ and $E_1$ have different norms, one can use the 
following upper estimate for $d_\theta[\ell_2^n(E_0),\ell_2^n(E_1)]$ 
in terms of type~$2$ constants which is taken from \cite[3.5]{kouba}: Let
 $[F_0,F_1]$ be a finite-dimensional interpolation couple. Then
\begin{equation}
\label{koubaestimate}
d_\theta[F_0,F_1] \le \type{2}(F_0')^{1-\theta} \cdot 
\type{2}(F_1')^\theta.
\end{equation}
Note that the estimate given in \eqref{koubaestimate} is slightly different 
from that in Kouba's work; we refer the reader to  \cite{dm98}  for the 
details. \\[10pt]
Using the simple fact
 that $\type{2}(\ell_2^n(E_i'))=\type{2}(E_i')$  (see e.\,g.
 \cite[11.12]{djt}), \eqref{koubaestimate} gives  
$d_\theta[\ell_2^n(E_0),\ell_2^n(E_1)] \le \type{2}(E_0')^{1-\theta} \cdot 
\type{2}(E_1')^\theta$.
Furthermore, by the duality of type and cotype (see e.\,g. 
\cite[11.10]{djt}) and the interpolative nature of the type~2 constants
 (see e.\,g. \cite[(3.8)]{tj})
$
\cotype{2}([E_0,E_1]_\theta) \le \type{2}([E_0',E_1']_\theta) 
\le \type{2}(E_0')^{1-\theta} \cdot \type{2}(E_1')^\theta
$.
Altogether we arrive at
\begin{cor}
Let $X_0,X_1$ be $n$-dimensional lattices and 
$[E_0,E_1]$ a finite-dimensional interpolation couple. Then for 
$0<\theta<1$ 
\begin{equation}
\label{notequal}
d_\theta[X_0(E_0),X_1(E_1)] \le \sqrt{2} \cdot \concave{2}(X_0)^{1-\theta} 
\cdot \concave{2}(X_1)^\theta \cdot (\type{2}(E_0')^{1-\theta} \cdot 
\type{2}(E_1')^\theta)^2. 
\end{equation}
\end{cor} 
\section{The finite-dimensional case in general}
Our estimates for $d_\theta[X_0(E_0),X_1(E_1);Y_0(F_0),Y_1(F_1)]$ 
are as follows:
\begin{prop}
\label{doublekouba}
Let $X_0,X_1$ and $Y_0,Y_1$ be $n$-dimensional and $m$-dimensional 
 lattices, respectively,  and $[E_0,E_1]$, $[F_0,F_1]$ two arbitrary
 finite-dimensional interpolation couples. Then for \mbox{$0<\theta<1$}
\begin{equation}
\label{epsilon}
\begin{split}
d_\theta & [X_0(E_0),X_1(E_1);Y_0(F_0),Y_1(F_1)] \\
  &\le 16 \cdot 
[(\concave{2}(X_0) \cdot 
\concave{2}(Y_0))^{1-\theta} (\concave{2}(X_1) \cdot 
\concave{2}(Y_1))^\theta]^{5/2} 
 \cdot t_\theta[E_0,E_1] \cdot t_\theta[F_0,F_1],
\end{split}
\end{equation}
where, if $G$ represents either $E$ or $F$,
$$
t_\theta[G_0,G_1]:=
\begin{cases}
\cotype{2}(G)^{5/2} &\text{ if $G=G_0=G_1$,} \\
(\type{2}(G_0')^{1-\theta} \cdot \type{2}(G_1')^\theta)^{7/2} &\text{ else.}
\end{cases}
$$ 
\em
The proof is based on the following ``factorization lemma'' 
which will enable us 
to use the estimates from the Hilbert space case derived in \eqref{equal} and 
\eqref{notequal} in order to obtain estimates for the general case. 
As usual we denote by $\Gamma_2$ the Banach operator ideal of all operators
 $T$ which allow a factorization $T=RS$ through a Hilbert space, 
together with the norm
$ 
\gamma_2(T):=\inf \norm{R} \cdot \norm{S}.
$
\end{prop}   
\begin{lemma}
\label{splitting}
Let $[E_0,E_1]$ and $[F_0,F_1]$ be finite-dimensional interpolation couples.
 Then for $0<\theta<1$
$$
\norm{\Gamma_2([E_0,E_1]_\theta',[F_0,F_1]_\theta) \Id 
[\Gamma_2(E_0',F_0), \Gamma_2(E_1',F_1)]_\theta} \le 
d_\theta[E_0,E_1] \cdot d_\theta[F_0,F_1].
$$
\end{lemma}
\proof Let $T: [E_0,E_1]_\theta' \rightarrow 
[F_0,F_1]_\theta$ factorize as follows:
\begin{center}
\settriparms[1`1`1;500]
\Vtriangle[[E_0,E_1]_\theta'`[F_0,F_1]_\theta`\ell_2;T`R`S],
\end{center}
 and consider by bilinear interpolation the norm $1$ mapping
$$
U \times V \rightarrow W, \quad (u,v) \mapsto v \circ u',
$$
where
$$U:= [\LL(\ell_2,E_0), \LL(\ell_2,E_1)]_\theta, \qquad 
 V:=[\LL(\ell_2,F_0), \LL(\ell_2,F_1)]_\theta
$$ 
and
$$W:= [\Gamma_2(E_0',F_0), \Gamma_2(E_1',F_1)]_\theta.$$
Then
$$
\norm{T}_W = \norm{SR}_W \le \norm{R'}_U \cdot \norm{S}_V 
 \le d_\theta[E_0,E_1] \cdot d_\theta[F_0,F_1] \cdot \norm{R'} \cdot 
\norm{S}, 
$$
which clearly gives  
$\norm{T}_W \le d_\theta[E_0,E_1] \cdot d_\theta[F_0,F_1] \cdot \gamma_2(T)$.
 \qed
\par Another ingredient needed for the proof of 
Proposition~\ref{doublekouba}, is a simple estimate for the cotype~2~constant
 of vector-valued Banach function spaces. We omit its straightforward proof
 (which needs arguments already used  in the proof of 
Lemma~\ref{factorization}).
\begin{lemma}
\label{vcotype}
Let $X$ be a $2$-concave Banach function space and
 $E$ a Banach space of cotype~$2$. Then $X(E)$ has cotype~$2$, and 
$\cotype{2}(X(E)) \le \sqrt{2} \cdot \concave{2}(X) \cdot \cotype{2}(E)$. 
\end{lemma}
With this the proof of Proposition~\ref{doublekouba} is straightforward:
\\[10pt]
{\em Proof} of Proposition~\ref{doublekouba}. 
For the moment denote by $D_\theta$ the norm of the embedding
$$
\Gamma_2([X_0(E_0),X_1(E_1)]_\theta',
[Y_0(F_0),Y_1(F_1)]_\theta) \Id [\Gamma_2(X_0(E_0)', Y_0(F_0)),
 \Gamma_2(X_1(E_1)',Y_1(F_1))]_\theta
$$
and $d_\theta:=d_\theta[X_0(E_0),X_1(E_1);Y_0(F_0),Y_1(F_1)]$.
Using Pisier's Factorization Theorem (\cite[4.1]{pisier86} or 
\cite[31.4]{df}), the Calder\'{o}n formula \eqref{calderonformula}, 
Lemma~\ref{vcotype} and the interpolative nature of the $2$-concavity 
constants (see \eqref{concaveinterpol}) one has 
\begin{align*}
d_\theta &\le (2 \cdot \cotype{2}([X_0(E_0),X_1(E_1)]_\theta)
 \cdot \cotype{2}([Y_0(F_0),Y_1(F_1)]_\theta))^{3/2}  \cdot D_\theta\\
& = (2 \cdot \cotype{2}((X_0^{1-\theta} X_1^\theta)([E_0,E_1]_\theta))
 \cdot \cotype{2}((Y_0^{1-\theta} Y_1^\theta)([F_0,F_1]_\theta))
)^{3/2} \cdot D_\theta\\
& \le 8 \cdot (\concave{2}(X_0^{1-\theta} X_1^\theta) \cdot 
\concave{2}(Y_0^{1-\theta} Y_1^\theta) \cdot 
\cotype{2}([E_0,E_1]_\theta) \cdot 
\cotype{2}([F_0,F_1]_\theta))^{3/2} \cdot D_\theta \\
& 
\le 8 \cdot ((\concave{2}(X_0) \cdot 
\concave{2}(Y_0))^{1-\theta} \cdot (\concave{2}(X_1) \cdot 
\concave{2}(Y_1))^\theta \\
& \qquad \cdot \cotype{2}([E_0,E_1]_\theta) \cdot 
\cotype{2}([F_0,F_1]_\theta))^{3/2} \cdot D_\theta.
\end{align*}
Now the  estimates stated in the proposition follow from 
 Lemma~\ref{splitting} together with \eqref{equal} and \eqref{notequal}. 
\qed
\section{The proof of the theorem}
In order to prove the theorem we need some additional
notation. For a 
$\sigma$-finite measure space $(\Omega,\Sigma,\mu)$ let  $FIN_\chi(\mu)$
 be the set of all subspaces of
$S(\mu)$---the linear space of all
 strictly simple functions---which are 
generated by a finite sequence of characteristic functions of measurable,
 pairwise disjoint sets with finite non-zero measures, and with $S(\mu,E)$ 
we denote the linear space of all strictly simple functions  with values in 
a normed space $E$.
\par
Now let us start the {\em Proof} of the theorem: First observe that
 if we define 
$$
\mathcal{A} := \{ U(M) \, | \, U \in FIN_\chi(\mu), M \in FIN(E_\Delta) \} 
$$
and
$$
 \mathcal{B}:=\{V(N) \, | \, V \in FIN_\chi(\nu), N \in 
FIN(F_\Delta)\},
$$ 
then $([X_0(E_0),X_1(E_1)], S(\mu,E_\Delta), \mathcal{A})$ and 
$ ([Y_0(F_0),Y_1(F_1)],S(\nu,F_\Delta),\mathcal{B})$ are cofinal interpolation 
triples whenever $X_0,X_1$ and $Y_0,Y_1$ have non-trivial concavity. 
Indeed, these assumptions together
 with \cite[1.a.5]{lt} and \cite[1.a.7]{lt} imply that $X_0$ and 
$X_1$ are $\sigma$-order continuous, and by \cite[p.~211]{kps} it 
follows that $S(\mu,E_\Delta)$ is dense in $X_0(E_0)$ and $X_1(E_1)$; 
 obviously each $G \in FIN(S(\mu,E_\Delta))$ is contained in some 
$U(M)$ with $U \in FIN_\chi(\mu)$ and $M \in FIN(E_\Delta)$. Moreover,
if $U$ is generated by measurable, pairwise disjoint sets
 $A_1, \ldots,A_n$ with finite non-zero measures, then $\chi_{A_1}, \ldots, 
\chi_{A_n}$ form a $1$-unconditional basis for $U$, hence $U$ is a 
finite-dimensional lattice which is order isometric to  $\R^n$ endowed with 
a lattice norm under the canonical order. 
\\[10pt]
 This now puts us in the position to apply the 
Approximation Lemma~\ref{epscor} together with Proposition~\ref{doublekouba}.
 For $U \in FIN_\chi(\mu)$, 
$V \in FIN_\chi(\nu)$, $M\in FIN(E_\Delta)$ and $N \in FIN(F_\Delta)$
\begin{align*}
 d_\theta  [U_0(M_0),U_1(M_1);& V_0(N_0),  V_1(N_1)] 
\\ &\le 16 \cdot 
[(\concave{2}(U_0) \cdot 
\concave{2}(V_0))^{1-\theta} (\concave{2}(U_1) \cdot 
\concave{2}(V_1))^\theta]^{5/2} \\
& \quad \cdot t_\theta[M_0,M_1] \cdot t_\theta[N_0,N_1] \\
& \le 16 \cdot 
[(\concave{2}(X_0) \cdot 
\concave{2}(Y_0))^{1-\theta} (\concave{2}(X_1) \cdot 
\concave{2}(Y_1))^\theta]^{5/2} \\
& \quad \cdot t_\theta[E_0,E_1] \cdot t_\theta[F_0,F_1],
\end{align*}
where the latter inequality follows from the fact that $\concave{2}$ 
respects sublattices, $\cotype{2}$  subspaces and $\type{2}$  quotients.
 \qed
%

\begin{thebibliography}{KPS82}

\bibitem[BL78]{BL}
J.~Bergh and J.~L\"ofstr\"om, \emph{Interpolation spaces}, Springer-Verlag,
  1978.

\bibitem[Cal64]{calderon}
A.~P. Calder\'{o}n, \emph{Intermediate spaces and interpolation, the complex
  method}, Studia Math. \textbf{24} (1964), 113--190.

\bibitem[Def99]{defant}
A.~Defant, \emph{Variants of the {M}aurey--{R}osenthal theorem for quasi 
{K}{\"o}the function spaces}, preprint (1999).

\bibitem[DF93]{df}
A.~Defant and K.~Floret, \emph{Tensor norms and operator ideals},
  North-Holland, 1993.

\bibitem[DM98]{dm98}
A.~Defant and C.~Michels, \emph{Bennett--{C}arl inequalities for symmetric
  {B}anach sequence spaces and unitary ideals}, submitted (1998).

\bibitem[DJT95]{djt}
J.~Diestel, H.~Jarchow, and A.~Tonge, \emph{Absolutely summing operators},
  Cambridge Studies in Advanced Mathematics 43, 1995.

\bibitem[Jun96]{junge}
M.~Junge, \emph{Factorization theory for spaces of operators}, Univ. Kiel,
  Habilitationsschrift, 1996. Currently available on the net under
\\ 
{\tt http://www-computerlabor.math.uni-kiel.de/\symbol{126}mjunge/preprints.html}

\bibitem[Kou91]{kouba}
O.~Kouba, \emph{On the interpolation of injective or projective tensor products
  of {B}anach spaces}, J. Funct. Anal. \textbf{96} (1991), 38--61.

\bibitem[KPS82]{kps}
S.~G. Krein, Ju.~I. Petunin, and E.~M. Semenov, \emph{Interpolation of linear
  operators}, Transl. Amer. Math. Soc. 54, 1982.

\bibitem[LT77]{lt77}
J.~Lindenstrauss and L.~Tzafriri, \emph{Classical {B}anach spaces {I}: Sequence
  spaces}, Springer-Verlag, 1977.

\bibitem[LT79]{lt}
J.~Lindenstrauss and L.~Tzafriri, \emph{Classical {B}anach spaces {II}:
  Function spaces}, Springer-Verlag, 1979.

\bibitem[Mau74]{maurey}
B.~Maurey, \emph{Th\'eor\`emes de factorisation pour les op\'erateurs
  lin\'eaires \`a valeurs dans les espaces {$L^p$}}, Ast\'erisque 11, 1974.

\bibitem[Pi86]{pisier86}
G.~Pisier, \emph{Factorization of linear operators and geometry of {B}anach
  spaces}, CBMS Regional Conf. Series 60, Amer. Math. Soc., 1986.

\bibitem[Pi90]{pisier90}
G.~Pisier, \emph{A remark on {$\Pi_2(\ell_p, \ell_p)$}}, Math. Nachr.
  \textbf{148} (1990), 243--245.

\bibitem[Pi96]{pisier96}
G.~Pisier, \emph{The operator {H}ilbert space {$OH$}, complex interpolation and
  tensor norms}, Mem. Amer. Math. Soc. 585, 1996.

\bibitem[TJ89]{tj}
N.~Tomczak-Jaegermann, \emph{Banach--{M}azur distances and finite-dimensional
  operator ideals}, Longman Scientific \& Technical, 1989.

\end{thebibliography}
%
\providecommand{\bysame}{\leavevmode\hbox to3em{\hrulefill}\thinspace}

Current address of both authors:
\\[5pt]
Fachbereich Mathematik \\
Carl von Ossietzky Universit\"at Oldenburg \\
Postfach 2503 \\
D-26111 Oldenburg \\
Germany 
\\[5pt]
\begin{tabular}{@{}ll}
E-mail: & {\tt defant}\symbol{64}{\tt mathematik.uni-oldenburg.de} \\
        & {\tt michels}\symbol{64}{\tt mathematik.uni-oldenburg.de}
\end{tabular}
\end{document}